\documentclass{article}
\pdfoutput=1

\usepackage{amsmath}
\usepackage[utf8]{inputenc}
\usepackage[english]{babel}

\usepackage{amssymb}
\usepackage{amsthm}
\usepackage{mathtools}
\usepackage{amsfonts}
\usepackage{enumerate}
\usepackage{enumitem}
\usepackage{multicol}
\usepackage{float}

\usepackage{tikz}
\usetikzlibrary{cd}

\allowdisplaybreaks

\newcommand{\ga}{\Gamma}

\theoremstyle{plain}
\newtheorem{theorem}{Theorem}[section]

\newtheorem{lemma}[theorem]{Lemma}
\newtheorem{cor}[theorem]{Corollary}
\newtheorem{remark}[theorem]{Remark}
\newtheorem{notation}[theorem]{Notation}

\newtheorem*{note}{Note}
\theoremstyle{definition}
\newtheorem{mydef}[theorem]{Definition}
\newtheorem{example}[theorem]{Example}

\author{}

\DeclareMathOperator{\Star}{St}
\DeclareMathOperator{\Link}{Lk}

\title{Geodesic Growth of some 3-dimensional RACGs}
\author{Yago Antol\'{i}n, Islam Foniqi }

\date{}

\begin{document}
	
	\maketitle
	\begin{abstract}
	We give explicit formulas for the geodesic growth series of a right-angled Coxeter group based on  a link-regular graph that does not contain 4-cliques.
    \end{abstract}
    \vspace{0.3cm}

\renewcommand{\thefootnote}{\fnsymbol{footnote}} 
\footnotetext{
\noindent\emph{MSC 2020 classification:} 05A15, 20F55.

\emph{Key words:} geodesic growth, Coxeter groups, right-angled Coxeter groups, rational formal power series, link-regular, clique.}     
\renewcommand{\thefootnote}{\arabic{footnote}} 
	
	\section{Introduction}
    Let $G$ be a group and $X$ a finite (monoid) generating set of $G$, that is, there exists a monoid morphism  $\pi\colon X^*\to G$ that is surjective.
Given a word $w\in X^*$, we denote by $\ell(w)$ its length.
A word $w\in X^*$ is a {\it geodesic} if $\ell(w)=\min_{u\in X^*}\{ \ell(u) : \pi(w) = \pi(u)\}$. The {\it geodesic growth series} associated to $(G,X)$ is the formal power series 
$$\mathcal{G}_{(G,X)}(z)= \sum_{i=0}^\infty \sharp \{w\in X^n : \text{ $w$ geodesic}\} z^n =\sum_{i=0}^n g_i z^i \in \mathbb{Z}[[z]].$$
One can similarly define the {\it standard growth series }associated to $(G,X)$ as the formal power series 
$$\mathcal{S}_{(G,X)}(z)= \sum_{i=0}^\infty \sharp \pi(\{w\in X^n : \text{ $w$ geodesic}\} ) z^n  =\sum_{i=0}^n s_i z^i\in \mathbb{Z}[[z]].$$
Both the standard and geodesic growth series encode  geometric information of the Cayley graph of $G$ (with respect to the generating set $X$). 
For instance, the celebrated theorem of Gromov on groups of polynomial growth says that the sequence of coefficients of $\{s_i\}_{i=0}^\infty$ is bounded by a polynomial on $i$ if and only if $G$ is virtually nilpotent. 
The geodesic growth has been less studied compared with standard growth, and moreover it is much more sensitive to the change of generating sets. In this paper, we focus on families of groups with some preferred generating sets.

{\it Right-angled Coxeter groups} (RACGs for short) is a family of groups described in terms of their defining presentation. Giving a simplicial graph $\Gamma$ with vertex set $V$ and edge set $E$, one associates to $\Gamma$ the right-angled Coxeter group $C_\Gamma$ defined by the following presentation:
$$C_\Gamma = \langle V \mid v^2=1 \, \forall v\in V,\; uv=vu \; \forall\, \{u,v\}\in E\rangle.$$
One calls $V$ the {\it standard generating set} of $C_\Gamma$.
We can also associate to $\Gamma$ the {\it right-angled Artin group} $A_\Gamma$ given by the presentation:
$$A_\Gamma = \langle V \mid uv=vu \; \forall\, \{u,v\}\in E\rangle.$$
One calls  $V\sqcup V^{-1}$  the {\it standard generating set } of $A_\Gamma$.
 
The languages of geodesics and shortlex representatives of a RACG with respect to its standard generating sets are regular \cite{BjornerBrenti,BrinkHowlett,LoefflerMeierWorthington}, and thus the corresponding standard and geodesic growth series are rational functions.
Concrete formulas for the  standard growth series of Coxeter groups, proved without the use of automata theory, can be found in \cite{Paris,Steinberg}.
Recently,  it was shown that  the growth rates of the geodesic and the standard  growth functions (i.e. $\lim \sqrt[n]{s_n}$ and $\lim \sqrt[n]{g_n}$) are either 1 or Perron numbers (see \cite{KolpTalam}).
 
Moreover, it is well understood how the geometry of the defining graph reflects on  the standard growth of $C_\Gamma$ with respect to $V$: it only depends on the cliques of $\Gamma$ of each size (see \cite[Proposition 17.4.2.]{DavisBook} or \cite{Paris,Steinberg}). 
For example, in the case of RACGs based on trees, this implies that the growth only depends on the number of vertices and edges of the tree.
 
Geodesic growth is still a very mysterious object compared to the standard growth and it is not clear which properties of the defining graph are reflected into the geodesic growth function. In \cite{CiobanuKolpakov}, Ciobanu and Kolpakov showed that there exist  infinitely many pairs of non-isomorphic RACGs based on trees with the same geodesic growth series with respect to the standard generators. These examples were based on co-spectral defining graphs, but then they gave infinitely many pairs of non-isomorphic RACGs with co-spectral defining graphs and different geodesic growth  series with respect to the standard generating set.
 
On the other hand, if the defining graph poses enough symmetry (the graph is link-regular), then the main theorem of \cite{antolin2013geodesic} states that the geodesic growth only depends on the number of cliques of each size and the isomorphism types of the links of the cliques. 
A simplicial graph is link-regular, if the number of elements of the link of a clique only depends on the size of the clique (See Definition \ref{link-regular}).

For example, if $\Gamma$ is a   totally disconnected graph with $n$ vertices, then it is link regular, and $C_\Gamma$ is a free product of cyclic groups of order 2. It is well-known that 
$$\mathcal{G}_{(C_\Gamma, V)}(z)-1 =\frac{nz}{1-(n-1)z}.$$
Moreover, in \cite{antolin2013geodesic} it is computed explicitly the geodesic growth series of right-angled Coxeter groups based on link-regular graphs with $n$ vertices, vertices of degree $l$, and without triangles. For such cases, one obtains the following formula for the  growth series.
$$
\mathcal{G}_{(C_\Gamma, V)}(z)-1=\dfrac{nz(1+(2-l)z)}{1+(-n-l+3)z+(-2n+2+nl)z^2}
.$$
Note that if $l=0$, one recovers the formula for a totally disconnected graph.

In this paper we continue to explore the geodesic growth of RACGs based on link-regular graphs. Our main result is to provide an explicit formula for the geodesic growth, if the graph does not contain 4-cliques.

\begin{theorem}\label{thm intro}
Let $\Gamma$ be a link-regular graph with $n$ vertices, $l$-regular and let $q$ be the link-number of an edge (which is the same for any edge), and without $4$-cliques.  Then,
$$\mathcal{G}_{(C_\Gamma, V)}(z)-1=\frac{nz(1 + (5 - l - q)z + (lq - 3l + 6)z^2)}
{1 + (6-n-l-q)z + (nl + lq + qn - 5n - 3l - q + 11)z^2  + (3nl + 6 - nlq - 6n )z^3}. $$
\end{theorem}
One can check that letting $q=0$, one obtains the previous formula for triangle-free link-regular graphs.
    
\section{Definitions and notation}
Let $\Gamma$  be a finite simplicial graph. For a vertex $a$ in $\Gamma$, we denote by $\Star(a)$ the set:
$$\Star(a) = \{b \in V\Gamma\mid \{a,b\} \in E\Gamma\} \cup \{a\}.$$
Let $\sigma \subseteq V\ga$ be such that the vertices of $\sigma$ span a complete subgraph of $\Gamma$, then $\sigma$ is called a \textit{clique}.
If $\sigma$ is a clique with $k$-vertices, then we call it a $k$-clique. 
Sometimes we refer to $3$-cliques and $4$-cliques by triangles and tetrahedrons respectively.

The {\it link of a clique $\sigma$}, denoted by $\Link(\sigma)$, is  the set of vertices in $V\Gamma\setminus \sigma$ that are connected with every vertex in $\sigma$. That is, 
 $$\Link(\sigma)=\{v \in V\Gamma \setminus \sigma : \{v\} \cup \sigma \text{ spans a clique}\}.$$
The {\it star of $\sigma$}, denoted by $\Star(\sigma)$, is the set of vertices in $\Gamma$ that are connected with every vertex in $\sigma$. That is, 
$$\Star(\sigma)=\{v \in V\Gamma : \{v\} \cup \sigma \text{ spans a clique}\}.$$
These sets satisfy $\sigma \cup \Link(\sigma)= \Star(\sigma)$.

\begin{mydef} \label{link-regular}
A graph $\Gamma$ is called \textit{link-regular} if for any clique $\sigma \in \Gamma$, $|\Link(\sigma)|$ depends on $|\sigma|$ and not on $\sigma$ itself, i.e. if $\sigma_1, \sigma_2$ are cliques with $|\sigma_1| = |\sigma_2|$ then $|\Link(\sigma_1)| = |\Link(\sigma_2)|$.
\end{mydef}

In this paper we will consider graphs which do not contain tetrahedrons. Under this condition, a graph is link-regular if there are numbers $l$ and $q$ such that the graph is $l$-regular (there are  $l$ edges meeting at any vertex), and  any edge is contained in $q$ triangles.

We recall the main theorem of \cite{antolin2013geodesic}. Remember that the $f$-polynomial associated to $\Gamma$ is the polynomial
$$f_\Gamma (z)= \sum_{n=0}^{|V|} \sharp \{ \Delta \subseteq \Gamma : \text{ $\Delta$ is an $n$-clique}\}z^n,$$
and essentially records the number of cliques of each size.
\begin{theorem}\label{main_antolin}
Let $C_\Gamma$ be the RACG based on a link-regular graph $\Gamma$. The geodesic growth of $G$ is fully determined by the $f$-polynomial of $\ga$ 
and the set of pairs $\{(|\sigma|, |\Link(\sigma)|) : \sigma \ \textrm{a clique in} \ \ga\}$.
\end{theorem}

\begin{remark}
As noted in \cite{athreya2014growth}, there is a relationship between the sizes of cliques and the coefficients of the $f$-polynomial. So one has that the geodesic growth of a RACG based on a link-regular graph $\Gamma$ is fully determined by the $f$-polynomial of $\ga$.
\end{remark}
        
\subsection{The double of a graph}
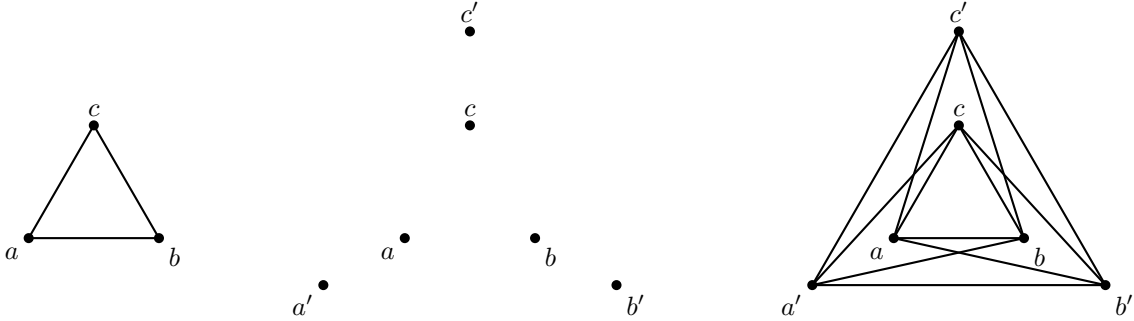
\begin{figure}[H]
\centering
\begin{tikzpicture}[thick]

\draw[fill=black] ({cos(90)}, {sin(90)}) circle (1.5pt) node[above] {$c$};

\draw[fill=black] ({cos(210)}, {sin(210)}) circle (1.5pt) node[below left] {$a$};

\draw[fill=black] ({cos(330)}, {sin(330)}) circle (1.5pt) node[below right] {$b$};

\draw[thick] ({cos(90)}, {sin(90)}) -- ({cos(210)}, {sin(210)});
\draw[thick] ({cos(210)}, {sin(210)}) -- ({cos(330)}, {sin(330)});
\draw[thick] ({cos(330)}, {sin(330)}) -- ({cos(90)}, {sin(90)});

\draw[fill=black] ({5 + cos(90)}, {sin(90)}) circle (1.5pt) node[above] {$c$};

\draw[fill=black] ({5 + cos(210)}, {sin(210)}) circle (1.5pt) node[below left] {$a$};

\draw[fill=black] ({5 + cos(330)}, {sin(330)}) circle (1.5pt) node[below right] {$b$};

\draw[fill=black] ({5 + 2.25*cos(90)}, {2.25*sin(90)}) circle (1.5pt) node[above] {$c'$};

\draw[fill=black] ({5 + 2.25*cos(210)}, {2.25*sin(210)}) circle (1.5pt) node[below left] {$a'$};

\draw[fill=black] ({5 + 2.25*cos(330)}, {2.25*sin(330)}) circle (1.5pt) node[below right] {$b'$};


\draw[fill=black] ({11.5 + cos(90)}, {sin(90)}) circle (1.5pt) node[above] {$c$};

\draw[fill=black] ({11.5 + cos(210)}, {sin(210)}) circle (1.5pt) node[below left] {$a$};

\draw[fill=black] ({11.5 + cos(330)}, {sin(330)}) circle (1.5pt) node[below right] {$b$};

\draw[thick] ({11.5 + cos(90)}, {sin(90)}) -- ({11.5 + cos(210)}, {sin(210)});
\draw[thick] ({11.5 + cos(210)}, {sin(210)}) -- ({11.5 + cos(330)}, {sin(330)});
\draw[thick] ({11.5 + cos(330)}, {sin(330)}) -- ({11.5 + cos(90)}, {sin(90)});

\draw[fill=black] ({11.5 + 2.25*cos(90)}, {2.25*sin(90)}) circle (1.5pt) node[above] {$c'$};

\draw[fill=black] ({11.5 + 2.25*cos(210)}, {2.25*sin(210)}) circle (1.5pt) node[below left] {$a'$};

\draw[fill=black] ({11.5 + 2.25*cos(330)}, {2.25*sin(330)}) circle (1.5pt) node[below right] {$b'$};

\draw[thick] ({11.5 + 2.25*cos(90)}, {2.25*sin(90)}) -- ({11.5 + 2.25*cos(210)}, {2.25*sin(210)});
\draw[thick] ({11.5 + 2.25*cos(210)}, {2.25*sin(210)}) -- ({11.5 + 2.25*cos(330)}, {2.25*sin(330)});
\draw[thick] ({11.5 + 2.25*cos(330)}, {2.25*sin(330)}) -- ({11.5 + 2.25*cos(90)}, {2.25*sin(90)});

\draw[thick] ({11.5 + 2.25*cos(90)}, {2.25*sin(90)}) -- ({11.5 + cos(210)}, {sin(210)});
\draw[thick] ({11.5 + 2.25*cos(210)}, {2.25*sin(210)}) -- ({11.5 + cos(330)}, {sin(330)});
\draw[thick] ({11.5 + 2.25*cos(330)}, {2.25*sin(330)}) -- ({11.5 + 2*cos(90)}, {sin(90)});

\draw[thick] ({11.5 + cos(90)}, {sin(90)}) -- ({11.5 + 2.25*cos(210)}, {2.25*sin(210)});
\draw[thick] ({11.5 + cos(210)}, {sin(210)}) -- ({11.5 + 2.25*cos(330)}, {2.25*sin(330)});
\draw[thick] ({11.5 + cos(330)}, {sin(330)}) -- ({11.5 + 2.25*cos(90)}, {2.25*sin(90)});

\end{tikzpicture}\caption{Construction of the double of a Graph: $\Gamma, V\Gamma^{[2]}$ and  $\Gamma^{[2]}$.}
\label{double}
\end{figure}

Given a graph $\Gamma = (V,E)$, the {\it double graph $\Gamma^{[2]}$} is defined as follows. Its vertex set is $V\Gamma^{[2]}=V\Gamma \sqcup V\Gamma$. Denote the vertices in the second copy as $\{a'\mid a\in V\}$.
For any edge $\{a,b\}$ in $\Gamma$ there are exactly four edges $\{a,b\}, \{a',b\}, \{a,b'\}, \{a',b'\}$ in $\Gamma^{[2]}$. See Figure \ref{double} for an example.

By definition, $|V\Gamma^{[2]}| = 2|V\Gamma|$, and $|E\Gamma^{[2]}| = 4|E\Gamma|$. One has also a projection $\rho:\Gamma^{[2]}\longrightarrow \Gamma$, which identifies naturally the two copies of the vertices of $\Gamma^{[2]}$. The map $\rho$ is $2$-to-$1$ on vertices and $4$-to-$1$ on edges. Moreover, by the construction of $\Gamma^{[2]}$ and the definition of $\rho$, we get $|\Link _{\Gamma^{[2]}}(v)| = 2|\Link_{\Gamma}(\rho(v))|$, i.e. $\rho$ is $2$-to-$1$ on links of vertices.

\begin{lemma} If $\Gamma$ is link-regular without tetrahedrons, then so is $\Gamma^{[2]}$.
\end{lemma}
\begin{proof}
We already have $|\Link _{\Gamma^{[2]}}(v)| = 2|\Link_{\Gamma}(\rho(v))|$.

For any edge $e = \{u,v\}$ in $\Gamma^{[2]}$ consider the edge $\rho(e)$ in $\Gamma$, and all the triangles $(\rho(u),\rho(v),c)$ over it. Triangles over $e = \{u,v\}$ are $(u,v,w)$ where $w\in \rho^{-1}(c)$. This means that $|\Link _{\Gamma^{[2]}}(e)| = 2|\Link _{\Gamma}(\rho(e))|$.

If there was a tetrahedron $\{x_1,x_2,x_3,x_4\}$ in $\Gamma^{[2]}$, then $\{\rho(x_1), \rho(x_2), \rho(x_3), \rho(x_4)\}$ would be a tetrahedron in $\Gamma$. Indeed $\rho(x_1), \rho(x_2), \rho(x_3), \rho(x_4)$ are all different since $\{x,x'\}$ cannot be an edge in $\Gamma^{[2]}$, and all the edges of the tetrahedron in $\{x_1,x_2,x_3,x_4\}$ induce edges for a tetrahedron over $\rho(x_1), \rho(x_2), \rho(x_3), \rho(x_4)$. Since there are no tetrahedrons in $\Gamma$, we conclude that there are no tetrahedrons in $\Gamma^{[2]}$.
\end{proof}

\begin{remark}
An important application of the double construction is provided in \cite[Lemma 2]{droms1993cayley}: one has that the Cayley graph of the RAAG based on $\ga$ is isomorphic as an undirected graph to the Cayley graph of the RACG based on $\ga^{[2]}$.
\end{remark}
Using the remark above we get the following:

\begin{cor}\label{raag-racg}
Let $G = A_\Gamma$ be a $RAAG$ based on  a link-regular graph $\Gamma$. The geodesic growth of $G$ is equal to the geodesic growth of $C_{\ga^{[2]}}$.
\end{cor}
        
\subsection{Geodesics in RACGs}
In this section we give characterizations of geodesics in RACGs.
Using Theorem 3.9 in \cite{wreo236}, and the characterization of a reduced sequence for graph products, we get the following result for RACGs.

\begin{theorem}\label{thm: Green geodesics}
Let $C_\Gamma$ be a right-angled Coxeter group based on $\Gamma$ and $V=V\Gamma$ the standard generating set.  Let  $w= (s_1,\ldots, s_n)$ be a word over $V$. 
Then $w$ is not a geodesic if and only if there are indices $1\leq i<j \leq n$ such that $s_i=s_j$ and $[s_i,s_k]=1$ for all $k$ satisfying $i<k<j$.
\end{theorem}

\begin{notation}\label{notation E}
Let $C_\Gamma$ be a RACG associated to $\ga$, with generating set $V = V\ga$.
If $w\in V^*$ is a word, we denote by $E_w$ the set of geodesics ending in $w$, and by $E_w(z)$ the generating growth series of $E_w$. That is
$$E_w(z)= \sum_{n=0}^\infty \sharp (E_w \cap V^n) z^n.$$
\end{notation}
\begin{theorem}\label{geodesics_in_racgs}
With the above notation.
Let $w_1,w_2$ be geodesic words over $V$ and $x\in V$, such that $w_1w_2x$ is also a geodesic. Assume $bx=xb$ in $C_\Gamma$, for all letters $b$ in $w_2$, and $x$ does not commute with any letter $a$ in $w_1$ (so, $w_1$ is not the empty word in particular). Then:
$$E_{w_1w_2x} = E_{w_1w_2}x, \text{ and also } E_{w_1w_2x}(z) = E_{w_1w_2}(z)\cdot z.$$
\end{theorem}

\begin{proof}
Obviously, one has the inclusion $E_{w_1w_2x} \subseteq E_{w_1w_2}x$.
To show $E_{w_1w_2x} \supseteq E_{w_1w_2}x$, we take a geodesic word $w\in E_{w_1w_2}$.
We suppose that $wx$ is not geodesic and derive a contradiction.
As $w_1w_2x$ is geodesic, there exist a shortest suffix of $wx$, say $w_0w_1w_2x$  that is not geodesic. Clearly,
$w_0$ is non-empty, and we can write it as $yu$, with $y$ a letter and $u$ a word (maybe empty).
As $yuw_1w_2x$ is not geodesic, but every proper subword is, we get  by the Theorem \ref{thm: Green geodesics} that $x=y$ and $x$ commutes with every letter of $u,w_1$ and $w_2$ which is the desired contradiction.
\end{proof}

\begin{notation}
Let $A,B$ be two subsets of a finitely generated free monoid $V^*$. 
If there is a bijection $f\colon A\to B$ that is length preserving (i.e. $\ell(f(a))=\ell(a)$ for all $a\in A$) we write $A\equiv B$.

In particular, if $A\equiv B$ then the corresponding growth series $A(z)$ and $B(z)$ are equal.
\end{notation}
For example, with Notation \ref{notation E}, if $a,b\in V$ commute, then $E_{ab}\equiv E_{ba}$ and $E_{ab}(z)=E_{ba}(z)$.
    
\section{Main Theorem}
Throughout the rest of the paper, $\Gamma$ will be a link-regular finite simple graph (with $n$ vertices) which does not contain tetrahedrons.
We denote by  $l$  the number of edges meeting at any vertex, and $q$ the number of triangles containing a fixed edge.

\begin{notation}\label{notation E2}
Denote by $G$, the group $C_\Gamma$, defining the RACG associated to $\ga$, with generating set $V = V\ga$.
We will use $\Delta\Gamma$ to denote the set of $3$-cliques of $\Gamma$.
In the following theorem we use the notation 

$$
\mathcal{E}_v(z) =  \sum_{a \in V\Gamma}E_{a}(z),\qquad
\mathcal{E}_e(z) =  \sum_{\substack{(a,b) \in (V\Gamma)^2\\ \{a,b\}\in E\Gamma}}E_{ab}(z),\qquad
\mathcal{E}_\Delta (z) =  \sum_{\substack{(a,b,c) \in (V\Gamma)^3\\ \{a,b,c\}\in \Delta\Gamma}}E_{abc}(z).
$$
\end{notation}

\begin{theorem}\label{thm1} 
Let $\Gamma$ be a link-regular graph with $n$ vertices, $l$-regular and let $q$ be the link-number of an edge (which is the same for any edge), and without $4$-cliques. Let $G_{\Gamma}$ be the corresponding right-angled Coxeter group, and $\mathcal{G}(z)$ the geodesic growth series of $G_\Gamma$ with respect to the standard generators.
Then, there exists polynomials $p_v, p_e, p_\Delta$ (given below) such that the following relations hold:
\begin{align}
    \mathcal{E}_v(z) = &
    \; \mathcal{G}(z) - 1 \label{E_v},\\
	\mathcal{E}_{e}(z) = &
	\; [\mathcal{G}(z) - 1](1-(n-l-1)z) - nz, \label{E_e}\\
	\mathcal{E}_{\Delta}(z) = &
	\; [\mathcal{G}(z) -1]p_{\Delta}(z, n, l, q) -nz +n(l-2q-2)z^2, \label{E_t}\\
	\sum_{(a,b,c,d) \in (V\Gamma)^4}E_{abcd}(z) = & \;
	[\mathcal{G}(z) -1] - nz - n(n-1)z^2 - [n(n-1)(n-2) + n(n-l-1)]z^3,\label{G_G}
\end{align}
and 
\begin{equation}
   \sum_{(a,b,c,d) \in (V\Gamma)^4}E_{abcd}(z) = [n + l + q - 6]z\cdot\mathcal{E}_{\Delta}(z) + p_{e}(n, l, q) z^2\cdot \mathcal{E}_{e}(z) + p_{v}(n, l, q) z^3\cdot \mathcal{E}_v(z). \label{E_4}
\end{equation}
Note than  one can find $\mathcal{G}(z)$ by substituting the  equations \eqref{E_v},\eqref{E_e},\eqref{E_t} and \eqref{G_G} into \eqref{E_4}.
Moreover, the polynomials $p_v,p_t,p_\Delta$ are given by:
\begin{align*}
    p_{\Delta}(z, n, l, q) = & 
    \; 1-(n+l-2q-3)z+(2(n-l-1)(l-q-1)-l(n-2l+q))z^2,\\
    p_{e}(n, l, q) = &
    \; n^2 +l^2 -2q^2 +  nl - 2nq  - 2lq   - 4n - 6l  + 10q + 7,\\
    p_{v}(n, l, q) = &
    \; (n-l-1)^3+2l(n-2l+q)(n-q-2)+lq(n-3l+3q).
\end{align*}
\end{theorem}

\begin{remark}
Equation \eqref{G_G} is obtained by subtracting from $\mathcal{G}(z)$ the generating growth series of geodesics of length at most $3$.
\end{remark}
\begin{note}
Theorem \ref{thm1} provides a way to calculate $\mathcal{G}(z)$ using a system of linear equations. The coefficients of the system are polynomials on $z$ (and $n,l,q$) of degree at most $3$. 
Theorem \ref{thm intro} was obtained by solving this linear system of equations with the help of Sage.

In the following example we calculate $\mathcal{G}(z)$ on a particular family. The general case can be computed similarly.
\end{note}
    
\begin{example}
Let's compute the geodesic growth of an infinite family $\Gamma_m$, defined inductively by 
$$\Gamma_0 = (V\Gamma_0= \{a,b,c\}, E\Gamma_0 = \{\{a,b\},\{a,c\},\{b,c\}\}) = \text{ triangle},$$ $$
\Gamma_{m+1} = \Gamma_m^{[2]} = \text{ the double of $\Gamma_m$.}$$

In terms of $(n,l,q)$ we have $(n_0,l_0,q_0) = (3,2,1)$, and therefore, by the double construction: 
$$(n_m,l_m,q_m) = (3 \cdot 2^m,2 \cdot 2^m,2^m).$$
\end{example}

Taking $2^m = k$ and substituting $3k,2k,k$ on the polynomials of Theorem \ref{thm1} for $n,l,q$ respectively, we find:
\begin{align*}
    p_{\Delta}(z, 3k, 2k, k) = &
    \; 2 (k - 1)^2 z^2 - 3( k - 1) z + 1, \\
    p_{e}(3k, 2k, k) = &
    \; 7 (k - 1)^2, \\
    p_{v}(3k, 2k, k) = &
    \; (k - 1)^3.
\end{align*}

Now substitute $3k, 2k, k$ on the other equations of Theorem \ref{thm1} for $n,l,q$ respectively. Use also the polynomials above, and we can find $\mathcal{G}(z)$ as a solution of the following system:
\begin{align*}
    \mathcal{E}_v(z) & = \mathcal{G}(z) - 1,\\
	\mathcal{E}_e(z) & = (\mathcal{G}(z) - 1)[1-(k-1)z] - 3kz,\\
	\mathcal{E}_\Delta (z) & = (\mathcal{G}(z) -1)[2 (k - 1)^2 z^2 - 3( k - 1) z + 1] -3kz - 6kz^2, \\
	\sum_{(a,b,c,d) \in (V\Gamma)^4}E_{abcd}(z) & =
	6(k - 1)z\mathcal{E}_\Delta (z) + [7 (k - 1)^2] z^2\mathcal{E}_e(z) + [(k - 1)^3] z^3\mathcal{E}_v(z),\\
	\mathcal{G}(z) & = 
	 1 + 3kz + 3k(3k-1)z^2 +[3k(3k-1)(3k-2) + 3k(k-1)]z^3 \\
	& \; \; \; \; + \sum_{(a,b,c,d) \in (V\Gamma)^4}E_{abcd}(z).
\end{align*}

Finally, solving for $\mathcal{G}(z)$, we find:

$$\mathcal{G}(z) = -\frac{6 z^3  + (2 k^2 - 7 k + 11) z^2 - 3 (k - 2) z + 1}{(z(k-1)-1)(2z(k-1)-1)(3z(k-1)-1)},$$

which agrees with the formula provided in Theorem \ref{thm intro}, for $(n,l,q) = (3k,2k,k)$.
    
Now using Theorem \eqref{thm intro}, and Corollary \eqref{raag-racg} we get the following:

\begin{cor}\label{cor1}
Let $\Gamma$ be a graph as in the hypothesis of Theorem \eqref{thm intro}. One can find the geodesic growth series $\mathcal{A}(z)$ for the  right-angled Artin group based on $\Gamma$ with respect to the generating set $V\Gamma \cup V\Gamma^{-1}$ by substituting $2n, 2l, 2q$ for $n,q,l$ respectively in the formula of Theorem \eqref{thm intro} and we get:
$$\mathcal{A}(z)-1=\frac{2nz[1 + (5 - 2l - 2q)z + (4lq - 6l + 6)z^2)]}
{1 + (6 - 2n - 2l - 2q)z + (4nl + 4lq + 4qn - 10n - 6l - 2q + 11)z^2 + (12nl + 6 - 8nlq - 12n)z^3}. $$
\end{cor}
    
\section{Proof of the main theorem}

Throughout this section $\Gamma$ is a link-regular graph without tetrahedrons. The graph $\Gamma$ has  $n$ vertices, the link of each vertex has $l$ vertices, and the link of each edge has $q$ vertices. 
Let $G = C_{\Gamma}$  be the associated RACG.

Note that there is 1 geodesic word of length 0, $n$ geodesic words of length $1$, $n(n-1)$ geodesics of length $2$. 
A word of length 3 is geodesic in $G$ if all its 3 letters are different or if it is of the form $aba$ with $b\notin \Star(a)$.
Thus there are $n(n-1)(n-2)+n(n-l-1)$ geodesic words of length $3$.

With Notation \ref{notation E}, one can write the geodesic growth series $\mathcal{G}(z)$ in any of the following forms:
\begin{align}
	\mathcal{G}(z) =&\; 1 + \sum_{a \in V\Gamma}E_a(z), 	\label{g_1} \\
	\mathcal{G}(z) =&\; 1 + nz + \sum_{(a,b) \in (V\Gamma)^2}E_{ab}(z), \label{g_2} \\
	\mathcal{G}(z) =&\; 1 + nz + n(n-1)z^2 + \sum_{(a,b,c) \in (V\Gamma)^3}E_{abc}(z), \label{g_3} \\
	\mathcal{G}(z) =&\; 1 + nz + n(n-1)z^2 +[n(n-1)(n-2) + n(n-l-1)]z^3 +\sum_{(a,b,c,d) \in (V\Gamma)^4}E_{abcd}(z). \label{g_4}
\end{align}
We get  
\eqref{E_v} and \eqref{G_G} of Theorem \ref{thm1} from equations \eqref{g_1} and \eqref{g_4}, respectively.
    
\begin{figure}[H]
\centering
\begin{tikzpicture}

\draw[fill=black] (0,0) circle (2pt);
\draw[fill=black] (2,0) circle (2pt);
\node at (0,0.3) {a};
\node at (2,0.3) {b};
\draw[thick] (0,0) -- (2,0);
\node at (1,-0.5) {(I)};

\draw[fill=black] (4,0) circle (2pt);
\draw[fill=black] (6,0) circle (2pt);
\node at (4,0.3) {a};
\node at (6,0.3) {b};

\node at (5,-0.5) {(II)};
\end{tikzpicture}

\caption{Configurations of 2 generators.} \label{fig:2v}
\end{figure}
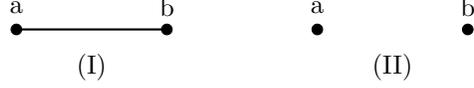

We will derive  \eqref{E_e} of the Theorem \ref{thm1} from \eqref{g_2} by expanding $\sum_{(a,b) \in (V\Gamma)^2}E_{ab}(z)$. Given a word $ab\in V^*$, we distinguish three cases: $a=b$,  $b\in \Link(a)$ and when $b\notin \Star(a)$.
We can describe these cases geometrically as in the Figure \ref{fig:2v} (omitting the case $a=b$).

The case $a=b$ is impossible, since no geodesic ends with $aa$.
In the case when $b\not\in \Star(a)$ we can write $E_{ab} = E_a \cdot b$. Hence $E_{ab}(z) = E_a(z) \cdot z$ and we have $n-l-1$ choices for $b$.

\begin{align}
\sum_{(a,b) \in (V\Gamma)^2}E_{ab}(z) & = \sum_{a\in V\Gamma} \left(\sum_{b\in \Link(a)}E_{ab}(z) + \sum_{b\not\in \Star(a)}E_{ab}(z)\right) \nonumber\\
& = \sum_{a\in V\Gamma} \left(\sum_{b\in \Link(a)}E_{ab}(z) \right) +  \sum_{a\in V \Gamma}(n-l-1)zE_{a}(z) \nonumber\\[5pt]
& = \mathcal{E}_e(z) + (n-l-1)z\mathcal{E}_v(z)\label{Eab}.
\end{align}

So, from equations \eqref{g_1}, \eqref{g_2}, and \eqref{Eab} we get:

\begin{equation*}
    \mathcal{E}_e(z) = (1-(n-l-1)z)[\mathcal{G}(z) - 1] - nz
\end{equation*}

which appears in the main theorem as the equation \eqref{E_e}.
    
We can work similarly to get \eqref{E_t} of the main theorem from the equation \eqref{g_3}. Consider the word $w = abc$. Since we will consider the geodesics that end in $w$, one needs $w$ itself to be a geodesic, and this implies that $a\neq b$, and $b \neq c$. The generators of the geodesic $abc$, lie in one of the following  disjoint cases:

\vspace{0.25cm}
\hspace{3cm}
\begin{minipage}{0.35\linewidth}
\begin{itemize}
    \item[(I)] $\{a,b\} \subseteq \Star(c)$
	\begin{itemize}
		\item[(I.1)] $a \in \Link(b)$
		\item[(I.2)] $a \not\in \Link(b)$
	\end{itemize}
\end{itemize}
\end{minipage}
\begin{minipage}{0.35\linewidth}
\begin{itemize}
    \item[(II)] $\{a,b\} \not \subseteq \Star(c)$
	\begin{itemize}
		\item[(II.1)] $a \in \Link(b)$
		\item[(II.2)] $a \not\in \Link(b)$
	\end{itemize}
\end{itemize}
\end{minipage}
\vspace{0.25cm}

\begin{figure}[H]
\centering
\begin{tikzpicture}[scale = 0.75]

\draw[fill=black] (0,0) circle (2pt);
\draw[fill=black] (2,0) circle (2pt);
\draw[fill=black] (0,2) circle (2pt);
\node at (0,2.3) {a};
\node at (2,-0.3) {b};
\node at (0,-0.3) {c};
\draw[thick] (0,0) -- (2,0) -- (0,2) -- (0,0);
\node at (1,-1) {(I.1)};

\draw[fill=black] (3,0) circle (2pt);
\draw[fill=black] (5,0) circle (2pt);
\draw[fill=black] (3,2) circle (2pt);
\node at (3,2.3) {a};
\node at (5,-0.3) {b};
\node at (3,-0.3) {c};
\draw[thick] (3,2) -- (3,0) -- (5,0);
\node at (4,-1) {(I.2)};

\draw[fill=black] (6,0) circle (2pt);
\draw[fill=black] (8,0) circle (2pt);
\draw[fill=black] (6,2) circle (2pt);
\node at (6,2.3) {a};
\node at (8,-0.3) {b};
\node at (6,-0.3) {c};
\draw[thick] (6,2) -- (8,0);
\draw[dotted][thick] (6,2) -- (6,0) -- (8,0);
\node at (7,-1) {(II.1)};

\draw[fill=black] (9,0) circle (2pt);
\draw[fill=black] (11,0) circle (2pt);
\draw[fill=black] (9,2) circle (2pt);
\node at (9,2.3) {a};
\node at (11,-0.3) {b};
\node at (9,-0.3) {c};
\draw[dotted][thick](9,2) -- (9,0) -- (11,0);
\node at (10,-1) {(II.2)};
\end{tikzpicture}
\caption{Configurations of 3 generators when $a\neq b \neq c$. Dashed edges might or might not appear in the configuration. Two vertices connected by a dashed edge might be the same vertex in $\Gamma$ (e.g. $a=c$ in (II.2)). No two dashed edges can be edges of the configuration simultaneously.} \label{fig:3v}
\end{figure}
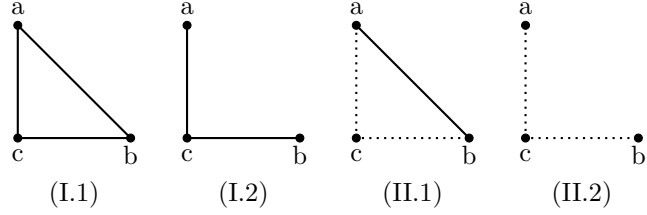

We can express these cases using the configurations of generators as in Figure \ref{fig:3v}. In the first three cases, the generators $a,b,c$ are all distinct, as explained in the respective cases. The generators in cases (I.1), and (I.2) appear in the defining graph $\Gamma$ exactly as they appear in Figure \ref{fig:3v}. 
In the case (II.2) one can have $a=c$ thought, and the way that $a,b,c$ appear in the defining graph $\Gamma$, in cases (II.1) and (II.2), depends on a subcase study. 


We can express $\sum_{(a,b,c) \in (V\Gamma)^3}E_{abc}(z)$ as a sum over the $4$ disjoint subcases given above. For convenience, we will write $\sum_X$, where $X$ is a case, to denote the summation over all triples $(a,b,c)\in (V\Gamma)^3$ satisfying the hypothesis of case $X$.
In (I.1) we also have $a\neq c$, because otherwise we would get $abc = aba = a^2b$ which would not be a geodesic. So, in this case the generators $a,b,c$ form a triangle and we get $\sum_{(\text{I.1})} E_{abc}(z) = \mathcal{E}_\Delta (z)$.
In (I.2) one gets $a\neq c$ as well, as $a = c$ would imply $\{a,b\} \subseteq \Star(a)$ and hence $a\in \Star(b)$ which cannot happen since $a\neq b$ and $a\not\in \Link(b)$. By Theorem \ref{geodesics_in_racgs} we have $E_{abc} \equiv E_{acb} = E_{ac}\cdot b$. 
Now we get $E_{abc}(z) = E_{acb}(z) = E_{ac}(z)\cdot z$, for a fixed $b$. Starting by fixing the edge $e=\{a,c\}$ we get $l-1-q$ choices for $b$, so in this case we have 
\begin{align*}
\sum_{\text{(I.2)}}E_{abc}(z)&= \sum_{a\in V}\sum_{c\in \Link(a)}\sum_{\substack{b\in \Link(c)\\b\notin \Link(a)}}E_{abc}(z) \\ &=
\sum_{a\in V}\sum_{c\in \Link(a)}(l-1-q)zE_{ac}(z) = (l-1-q)z\mathcal{E}_e(z).
\end{align*}
In (II.1) once again $a\neq c$, as $a = c$ would imply $\{a,b\} \not \subseteq \Star(a)$ and hence $a\not \in \Star(b)$ which cannot happen since $a\in \Link(b)$. We count by first fixing the edge $\{a,b\}$, and then letting $c$  be any vertex different to $a$ and $b$ that does not form a triangle with $\{a,b\}$. One has $n-2 -q$ choices for $c$. As both of $\{a,c\}, \{b,c\}$ cannot be edges, using Theorem \ref{geodesics_in_racgs} and considering all the subcases, we get the formula $E_{abc} = E_{ab}\cdot c$. Now, arguing as in (I.2) we get
$$\sum_{\text{(II.1)}}E_{abc}(z)=(n-q-2)z\mathcal{E}_e(z).$$
In (II.2) we count by first fixing the vertex $a$ and then considering the choices when $c=a, c\in \Link(a)$, and $c\not\in \Star(a)$.
In this case, $b\notin\Star(a)\cup \{c\}$ and moreover, when $c\in \Link(a)$ then $b$ is not linked to $c$. 

For $c=a$ we have $ (n-l-1)$ possible choices for $b$; for $c\in \Link(a)$  we have $l\cdot (n-2l+q)$ possible choices for $c,b$; finally for $c\notin \Star(a)$ we get  $(n-l-1)(n-l-2)$ possible choices for $c,b$.
Using Theorem \ref{geodesics_in_racgs} we get the formula $E_{abc} = E_{ab}\cdot c = E_{a}\cdot b \cdot c = E_{a}\cdot bc$, so 
\begin{align*}
\sum_{\text{(II.2)}}E_{abc}(z)&= 
\sum_{a\in V}\left( 
\sum_{c=a}
\sum_{b\notin \Link(a)} E_{abc}(z)
+\sum_{c\in \Link(a)} \sum_{\substack{b\notin \Link(a)\\b\notin \Link(c)}}E_{abc}(z) 
+\sum_{c\notin \Star(a)} \sum_{\substack{b\notin \Link(a)\\ b\neq c}} E_{abc}(z)\right)\\
&=\sum_{a\in V}\left((n-l-1)z^2E_a(z)+ l(n-2l+q) z^2E_a(z)+ (n-l-1)(n-l-2)z^2E_a(z)\right)\\
& =[(n-l-1)^2 + l(n-2l+q)]z^2\mathcal{E}_v(z).
\end{align*}
Now summing everything up we get:
\begin{equation}\label{Eabc} 
    \sum_{(a,b,c) \in (V\Gamma)^3}E_{abc}(z) = \mathcal{E}_{\Delta}(z) + (n+l-2q-3)z\mathcal{E}_e(z) +  [(n-l-1)^2 + l(n-2l+q)]z^2\mathcal{E}_v(z)
\end{equation}
Substituting \eqref{E_v},\eqref{E_e},\eqref{g_3} into \eqref{Eabc}, we get
\begin{align*}
 	\mathcal{G}(z) -( 1 + nz + n(n-1)z^2 ) & =  \mathcal{E}_{\Delta}(z) \\
 	& + (n+l-2q-3)\cdot z\cdot  \left((1-(n-l-1)z)[\mathcal{G}(z) - 1] - nz\right) \\
 	& +  [(n-l-1)^2 + l(n-2l+q)]\cdot z^2\cdot  (\mathcal{G}(z) -1)
\end{align*}
And one gets a formula for $\mathcal{E}_{\Delta}(z)$:
\begin{align*}
    \mathcal{E}_{\Delta}(z) = &\;(\mathcal{G}(z) - 1)[1-(n+l-2q-3)z+(2(n-l-1)(l-q-1)-l(n-2l+q))z^2]\\
	& -nz +n(l-2q-2)z^2
\end{align*}
which appears in the main theorem as Equation \eqref{E_t}.
    
To finish the proof, we need to show that \eqref{E_4} holds. As in the previous cases, we proceed to rewrite $\sum_{(a,b,c,d) \in (V\Gamma)^4}E_{abcd}(z)$ depending on different cases for the word $abcd$. Since we consider the geodesics that end in $abcd$, we want $abcd$ to be a geodesic itself, and this implies that $a\neq b$, $b \neq c$, and $c \neq d$.

We distinguish the following disjoint cases:

\vspace{0.25cm}
\hspace{3cm}
\begin{minipage}{0.35\linewidth}
\begin{itemize}
    \item[(I)] $\{a,b,c\} \subseteq \Star(d)$
		\begin{itemize}
			\item[(I.1)] $\{a,b\} \subseteq \Star(c)$
			\begin{itemize}
				\item[(I.1.1)] $a \in \Link(b)$
				\item[(I.1.2)] $a \not\in \Link(b)$		\end{itemize}
			\item[(I.2)] $\{a,b\} \not\subseteq \Star(c)$
			\begin{itemize}
				\item[(I.2.1)] $a \in \Link(b)$
				\item[(I.2.2)] $a \not\in \Link(b)$
			\end{itemize}
		\end{itemize}
    \end{itemize}
\end{minipage}
\begin{minipage}{0.35\linewidth}
\begin{itemize}
    \item[(II)] $\{a,b,c\} \not\subseteq \Star(d)$
		\begin{itemize}
			\item[(II.1)] $\{a,b\} \subseteq \Star(c)$
			\begin{itemize}
				\item[(II.1.1)] $a \in \Link(b)$
				\item[(II.1.2)] $a \not\in \Link(b)$
			\end{itemize}
			\item[(II.2)] $\{a,b\} \not \subseteq \Star(c)$
			\begin{itemize}
				\item[(II.2.1)] $a \in \Link(b)$
				\item[(II.2.2)] $a \not\in \Link(b)$
			\end{itemize}
		\end{itemize}
\end{itemize}
\end{minipage}
\vspace{0.25cm}

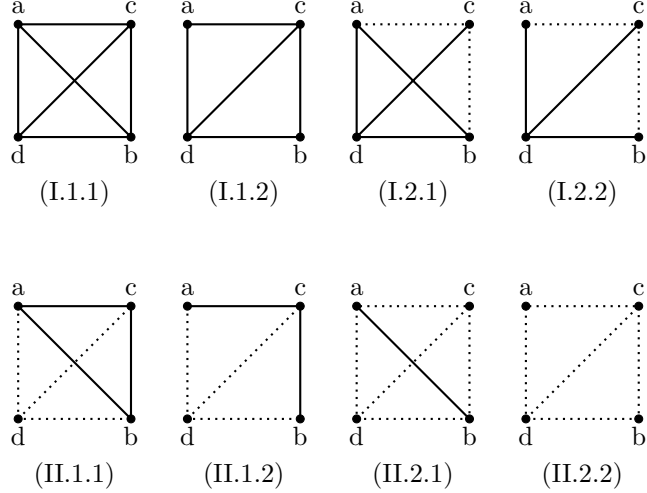
\begin{figure}[H]
\centering
\begin{tikzpicture}[scale=0.75]

\draw[fill=black] (0,0) circle (2pt);
\draw[fill=black] (2,0) circle (2pt);
\draw[fill=black] (0,2) circle (2pt);
\draw[fill=black] (2,2) circle (2pt);
\node at (0,2.3) {a};
\node at (2,-0.3) {b};
\node at (0,-0.3) {d};
\node at (2,2.3) {c};
\draw[thick] (0,0) -- (2,0) -- (2,2) -- (0,2) -- (0,0) -- (2,2);
\draw[thick] (0,2) -- (2,0);
\node at (1,-1) {(I.1.1)};

\draw[fill=black] (3,0) circle (2pt);
\draw[fill=black] (5,0) circle (2pt);
\draw[fill=black] (3,2) circle (2pt);
\draw[fill=black] (5,2) circle (2pt);
\node at (3,2.3) {a};
\node at (5,-0.3) {b};
\node at (3,-0.3) {d};
\node at (5,2.3) {c};
\draw[thick] (3,0) -- (5,0) -- (5,2) -- (3,2) -- (3,0) -- (5,2);
\node at (4,-1) {(I.1.2)};

\draw[fill=black] (6,0) circle (2pt);
\draw[fill=black] (8,0) circle (2pt);
\draw[fill=black] (6,2) circle (2pt);
\draw[fill=black] (8,2) circle (2pt);
\node at (6,2.3) {a};
\node at (8,-0.3) {b};
\node at (6,-0.3) {d};
\node at (8,2.3) {c};
\draw[thick] (6,0) -- (8,0) -- (6,2) -- (6,0) -- (8,2);
\draw[dotted][thick] (8,0) -- (8,2) -- (6,2);
\node at (7,-1) {(I.2.1)};

\draw[fill=black] (9,0) circle (2pt);
\draw[fill=black] (11,0) circle (2pt);
\draw[fill=black] (9,2) circle (2pt);
\draw[fill=black] (11,2) circle (2pt);
\node at (9,2.3) {a};
\node at (11,-0.3) {b};
\node at (9,-0.3) {d};
\node at (11,2.3) {c};
\draw[thick] (11,0) -- (9,0) -- (9,2);
\draw[thick] (9,0) -- (11,2);
\draw[dotted][thick] (11,0) -- (11,2) -- (9,2);
\node at (10,-1) {(I.2.2)};

\draw[fill=black] (0,-5) circle (2pt);
\draw[fill=black] (2,-5) circle (2pt);
\draw[fill=black] (0,-3) circle (2pt);
\draw[fill=black] (2,-3) circle (2pt);
\node at (0,-2.7) {a};
\node at (2,-5.3) {b};
\node at (0,-5.3) {d};
\node at (2,-2.7) {c};
\draw[thick] (2,-5) -- (2,-3) -- (0,-3) -- (2,-5);
\draw[dotted][thick] (0,-3) -- (0,-5) -- (2,-5);
\draw[dotted][thick] (0,-5) -- (2,-3);
\node at (1,-6) {(II.1.1)};

\draw[fill=black] (3,-5) circle (2pt);
\draw[fill=black] (5,-5) circle (2pt);
\draw[fill=black] (3,-3) circle (2pt);
\draw[fill=black] (5,-3) circle (2pt);
\node at (3,-2.7) {a};
\node at (5,-5.3) {b};
\node at (3,-5.3) {d};
\node at (5,-2.7) {c};
\draw[thick] (5,-5) -- (5,-3) -- (3,-3);
\draw[dotted][thick] (3,-3) -- (3,-5) -- (5,-5);
\draw[dotted][thick] (3,-5) -- (5,-3);
\node at (4,-6) {(II.1.2)};

\draw[fill=black] (6,-5) circle (2pt);
\draw[fill=black] (8,-5) circle (2pt);
\draw[fill=black] (6,-3) circle (2pt);
\draw[fill=black] (8,-3) circle (2pt);
\node at (6,-2.7) {a};
\node at (8,-5.3) {b};
\node at (6,-5.3) {d};
\node at (8,-2.7) {c};
\draw[dotted][thick] (6,-5) -- (8,-5) -- (8,-3) -- (6,-3) -- (6,-5) -- (8,-3);
\draw[thick] (6,-3) -- (8,-5);
\node at (7,-6) {(II.2.1)};

\draw[fill=black] (9,-5) circle (2pt);
\draw[fill=black] (11,-5) circle (2pt);
\draw[fill=black] (9,-3) circle (2pt);
\draw[fill=black] (11,-3) circle (2pt);
\node at (9,-2.7) {a};
\node at (11,-5.3) {b};
\node at (9,-5.3) {d};
\node at (11,-2.7) {c};
\draw[dotted][thick] (9,-5) -- (11,-5) -- (11,-3) -- (9,-3) -- (9,-5) -- (11,-3);
\node at (10,-6) {(II.2.2)};
\end{tikzpicture}
\caption{Configurations of 4  generators with $a \neq b \neq c \neq d$. Dashed edges represent pairs of vertices that might or might not be in the star of each other (including the possibility of being equal).  In cases (I.2) at least one dashed edge is not an edge in $\Gamma$ (i.e. the vertices are different and do not commute). In cases (II) at least one of the dashed edges incident to $d$ is not an edge in $\Gamma$. Moreover, in cases (II.2), 
at least one of $a$ or $b$ does not form an edge with $c$.}\label{fig:4v}
\end{figure}

We can express them geometrically as configurations of $4$ points in  Figure \ref{fig:4v}. Each individual figure will be considered in detail, as most of them represent a family of subcases, and not necessarily a truthful configuration of 4 generators in $\Gamma.$

As above, we will write $\sum_X$, where $X$ is one of the  cases above, to denote the summation over all quadruples $(a,b,c,d)\in (V\Gamma)^4$ satisfying the hypothesis of case $X$.
	
In (I.1.1) any two of the vertices $a,b,c,d$ would commute, so they would have to be pairwise distinct and hence form a tetrahedron; however in $\Gamma$ there are no $4$-cliques, so we obtain $\sum_{\text(I.1.1)}E_{abcd}(z)=0.$

In (I.1.2), except the pair $(a,b)$ any two other pairs of vertices among $a,b,c,d$ commute, so all of them have to be pairwise distinct, as otherwise $abcd$ would not be a geodesic. Fixing the triangle $\{a,c,d\}$,  we have $q-1$ choices for $b$. Also, for $a,b,c,d$ in this case, $E_{abcd} \equiv  E_{acdb} = E_{acd}\cdot b$ (by Theorem \ref{geodesics_in_racgs}, as $b$ does not commute with $a$), hence $E_{abcd}(z)=E_{acd}(z)\cdot z$. Therefore 

$$
\sum_{\text{(I.1.2)}}E_{abcd}(z) = (q-1)z\mathcal{E}_{\Delta}(z).
$$

In (I.2.1), the vertex $d$ commutes with any of the letters of the word $abcd$ so it should be distinct from any of $a,b,c$. The only equal pair could be $(a,c)$; but $a = c$ implies $\{a,b\}\not \subseteq \Star(a)$ which contradicts the assumption $a \in \Link(b)$. So, once again all the vertices are pairwise distinct. Fixing the triangle $\{a,b,d\}$, we have $l-2$ choices for $c$. Also, for $a,b,c,d$ in this case, we have $E_{abcd}\equiv E_{abdc} = E_{abd}\cdot c$, as $c$ does not commute with at least one of $a,b$. So, $E_{abcd}(z)=E_{abd}(z)\cdot z$, and now we obtain

$$
\sum_{\text{(I.2.1)}}E_{abcd}(z)
=(l-2)z\mathcal{E}_{\Delta}(z).$$

In (I.2.2), as in (I.2.1), the vertex $d$ commutes with any of the letters of the word $abcd$ so it should be distinct from any of $a,b,c$. Except the pair $(a,c)$, which could be equal, any other two vertices among $a,b,c,d$ are different. So $\{a,d\}$ forms an edge and we start by fixing it.
As, $b,c\in \Link(d)$, and as both $b$ and $c$ do not commute with at least one of $a, d$ we have $E_{abcd}\equiv E_{adbc} = E_{adb}\cdot c = E_{ad}\cdot b \cdot c$. Now one gets $E_{abcd}(z)=E_{ad}(z)\cdot z^2$.
We need to count the possibilities for  $b$ and $c$; consider the following disjoint subcases for $c$:
\begin{itemize}
	\item[(1)] $c=a:$ Here we have $l-1-q$ choices for $b$.
	\item[(2)] $c \in \Link(a):$ Here we have $q$ choices for $c$. Note that as $\{a,b\}\not\subseteq \Star(c)$, we obtain that $b\notin \Star(c)$. Therefore, $b$ is the link of $d$ but $b$ is not in the link of the edge $\{a,d\}$ and neither in the link of the edge $\{c,d\}$.
	Note that the $\Link(\{a,d\})\cap \Link(\{c,d\})$ is empty, as if not, $\Gamma$ will contain a tetrahedron. 
	Therefore, there are  $l-2q$ choices for $b$ in this subcase.
	\item[(3)] $c \not \in \Star(a):$ we have $l-1-q$ choices for $c$ and $l-2-q$ choices for $b$, since on top we have that $b\neq c$, as there is no geodesic of the form $accd$.
\end{itemize}
At the end one gets 
\begin{align*}
\sum_{\text{(I.2.2)}}E_{abcd}(z)&=
\sum_{a\in V}\sum_{d\in \Link(a)}
\left(
\sum_{c=a}E_{abcd}(z) 
+\sum_{c\in \Link(a)}E_{abcd}(z)
+\sum_{c\notin\Star(a)}E_{abcd}(z)
\right)  \\[5pt]
&=
\sum_{a\in V}\sum_{d\in \Link(a)}
\left(
(l-1-q)z^2E_{ad}(z) 
+q(l-2q)z^2 E_{ad}(z)
+(l-1-q)(l-2-q)z^2E_{ad}(z)\right)\\
&= [(l-1-q)^2 + q(l-2q)]z^2\mathcal{E}_e(z).
\end{align*}

We now consider the case (II).	

In (II.1.1), all the vertices should be distinct, indeed, any pair among $a,b,c$ commute, so $a\neq c$, and any of them is distinct from $d$ as one can permute the letters of the subword $abc$ in $abcd$. We count by fixing the triangle $abc$. 
We have $E_{abcd}=E_{abc}\cdot d$, as $d$ does not commute with at least one of $a,b,c$. One gets also $E_{abcd}(z)=E_{abc}(z)\cdot z$. 
Since $\Gamma$ does not have tetrahedrons, there is no condition on $d$ except that $d\not \in \{a,b,c\}$. We have $n-3$ choices for $d$, so we obtain 

$$\sum_{\text{(II.1.1)}}E_{abcd}(z)= (n-3)z\sum_{\Delta \in \Delta \Gamma}E_{\Delta}(z).$$ 
In (II.1.2), $a,b,c$ are distinct as for $a = c$ the get $a\in \Link(b)$. As $abc = acb$, the only equal pair could be $(a,d)$, which will be discussed below. Here we count by first fixing the edge $\{a,c\}$ and then considering the different cases for $d$: $d\in \{a,c\}$, $d\in \Link(\{a,c\})$ or $d\notin \Star(\{a,c\})$. 
As $c$ commutes with $b$, $b$ does not commute with at least one of $a,c$, and $d$ does not commute with at least one of $a,c,b$, we have 
$E_{abcd}\equiv E_{acbd} = E_{acb}\cdot d = E_{ac}\cdot bd$. Hence, also $E_{abcd}(z)=E_{ac}(z)\cdot z^2$.
We now  count the choices for $b,d$:
\begin{itemize}
	\item[(1)] $d\in \{a,c\}$: The case $d=c$  is impossible. The case  $d=a$, we have that $b\in \Link(c)\setminus\Star(a)$, and this gives $l-1-q$ choices for $b$. 
	\item[(2)] $d \in \Link(\{a,c\}):$  there are $q$ choices for $d$. In this case we have that $b\in \Link(c)$. Since $\{a,b,c\}\not\subseteq \Star(d)$ and $d\in \Link(\{a,c\})$, we get $b\notin \Link(d)$.
	Also from the hypothesis, we get $b\notin \Link(a)$.
	Therefore $b$ is in $\Link(c)\setminus (\Link(\{a,d\})\cup \Link(\{c,d\}))$.   Note that the links of two  edges in a triangle  are disjoint as $\Gamma$ has no tetrahedrons. Thus, there are $l-2q$ choices for $b$.
	\item[(3)] $d \not\in \Star(\{a,c\}):$ We subdivide this case into two subcases:
	\begin{itemize}
		\item[(3.1)] $d\in \Star(c):$ 
		In this case, $d$ is in $\Star(c)\setminus \Star(\{a,c\})$ and there are $(l+1)-(q+2)=l-q-1$ possibilities for $d$.
		
		Note that in this case $b\neq d$ since otherwise $abcd$ would not be a geodesic. Thus we have that $b$ is in $\Star(c)\setminus ( \Star(\{a,c\})\cup\{d\})$, and we have  $l-2-q$ possibilities for $b$. 
		\item[(3.2)] $d\not \in \Star(c):$ As $d\in V\setminus \Star(c)$,  we have $n-l-1$ choices for $d$. As $b\in \Star(c)\setminus \Star(\{a,c\})$, we have  $l-1-q$ choices for $b$.
	\end{itemize}
\end{itemize}
Ultimately, we obtain 
\begin{align*}
 \sum_{\text{(II.1.2)}}E_{abcd}(z) 
 & =
 \sum_{a\in V}\sum_{c\in \Link(a)}\left(
 \sum_{d\in \{a,c\}}E_{abcd}(z) 
+\sum_{d\in \Link(\{a,c\})}E_{abcd}(z)
+\sum_{d\notin\Star(\{a,c\})}E_{abcd}(z)
 \right)\\
 & =
 \sum_{a\in V}\sum_{c\in \Link(a)}\left(
1 \cdot (l-1-q) z^2 E_{ac}(z)+ 
q (l-2q) z^2 E_{ac}(z) + (l-q-1)(n-q-3)z^2 E_{ac}(z)\right)
\\
 &=[(l-1-q)(n-q-2)+q(l-2q)]z^2\mathcal{E}_e(z).
\end{align*}
In (II.2.1), we start by fixing the edge $\{a,b\}$ and then considering the different cases for $d$: $d\in \{a,b\}$, $d\in \Link(\{a,b\})$ or $d\notin \Star(\{a,b\})$. 
The vertices $a,b,c$ are all different, as $abc = bac$. Since $\{a,b,c\} \not\subseteq \Star(d)$ we get $E_{abcd}=E_{abc}\cdot d$. Moreover, $\{a,b\} \not\subseteq \Star(c)$, so  $E_{abc}=E_{ab}\cdot c$. Putting these together, we get  $E_{abcd}=E_{ab}\cdot cd$, and hence $E_{abcd}(z)=E_{ab}(z)\cdot z^2$. We  count the  choices for $c$ and $d$.
\begin{itemize}
	\item[(1)] $d\in\{a,b\}$: If $d=a$, we have $\{a,b\}\subseteq \Star(d)$ as $a\in \Link(b)$. Since $\{a,b,c\} \not\subseteq \Star(d)$ we get $c\not \in \Star(d)$. Thus $c$ can be any vertex outside of $\Star(d)$, which means $n-l-1$ choices for $c$. 
	If $d=b$, the discussion is analogous and we can take for $c$ any vertex outside $\Star(b) = \Star(d)$ and we have $n-l-1$ choices for $c$ yet again.
	\item[(2)] $d \in \Link(\{a,b\}):$ Here we have $q$ choices for $d$. Since $\{a,b,c\}\not\subseteq \Star(d)$, $c$ can not be in  $\Link(d)$. Since $\{a,b\}\not\subseteq \Star(c)$, $c$ is not in  $\Link(\{a,b\})$.  Since $a,b,d$ form a triangle, and we do not have tetrahedrons, these two links are disjoint. There are $n-q-l$ choices for $c$.
	\item[(3)] $d \not \in \Star(\{a,b\}):$ Here both $c$ and $d$ are not in $\Star(\{a,b\})$ and moreover $c\neq d$ to have $abcd$ a geodesic. We have $d\in V\setminus \Star(\{a,b\})$ that gives $n-q-2$ choices for $d$ and we have $c\in V\setminus (\Star(\{a,b\})\cup \{d\})$ that gives  $n-q-3$ for $c$.
\end{itemize}
Summing up, we obtain:
$$\sum_{\text{(II.2.1)}}E_{abcd}(z)=[2(n-l-1) + q(n-q-l) + (n-q-2)(n-q-3)]z^2\mathcal{E}_e(z).$$ 
In (II.2.2), we first fix $a$ and then we consider different cases for $d$: $d=a, d\in \Link(a)$ or $d\not\in \Star(a)$. Note that since $\{a,b,c\} \not\subseteq \Star(d)$ we have:
$E_{abcd}=E_{abc}\cdot d$. Similarly since $\{a,b\} \not\subseteq \Star(c)$ we have $E_{abc}=E_{ab}\cdot c$, and ultimately, since $a\not \in \Link(b)$ we obtain $E_{ab}=E_{a}\cdot b$. Putting everything together we get $E_{abcd}=E_{a}\cdot bcd$
and hence $E_{abcd}(z)=E_{a}(z)\cdot z^3$. We count the choices for $b,c,d$. 
\begin{itemize}
	\item[(1)] $d=a:$ Here we have $1$ choice for $d$. We split this case into  the following disjoint subcases:
	\begin{itemize}
	    \item[(1.1)] $c=a$: this is impossible, since $abaa$ is not a geodesic.
	    \item[(1.2)] $c\in \Link(a):$ Here we have $l$ choices for $c$. 
	    Since $\{a,b,c\} \not\subseteq \Star(d)$, $d = a$ and $c\in \Link(a) = \Link(d)$, we have $b\not \in \Star(a) = \Star(d)$. Further, since $\{a,b\} \not\subseteq \Star(c)$ and $a\in \Link(c)$ it must be that $b\not \in \Link(c)$.
	    In this case, $b$ can be any vertex of $V\setminus (\Star(a)\cup\Star(c))$. As $|\Star(a)\cap\Star(c)|=q+2$, we have 
	    $n-2(l+1)+(q+2)=n-2l+q$ possibilities for $b$. 
	    \item[(1.3)] $c\not \in \Star(a):$ Here $b\not\in \Star(a)\cup \{c\}$. So we have $n-l-1$ choices for $c$, and  $n-l-2$ for $b$.
	\end{itemize}
	Accounting for (1.1), (1.2) and (1.3), for a given vertex $a\in V$ we have:
	
	$$\sum_{\substack{c,b,d\in  \text{(II.2.2)}\\d=a}} E_{abcd}(z)= [l(n-2l+q)+(n-l-1)(n-l-2)] z^3\mathcal{E}_v(z).$$
	\item[(2)] $d\in \Link(a):$  Here we have $l$ choices for $d$.  We divide now the analysis into  the following disjoint subcases:
	\begin{itemize} 
	    \item [(2.1)] $c\in \{a,d\}$:  The case  $c=d$ is impossible. 
		In the case $c=a$, we have one choice for $c$ and  $b$ can be any vertex of $V\setminus (\Star(a)\cup\Star(d))$. As $|\Star(a)\cap\Star(d)|=q+2$, we have $n-2l+q$ possibilities for $b$.
		\item[(2.2)] $c\in \Link(\{a,d\}):$ which gives us $q$ choices for $c$. We have a triangle $\{a,c,d\}$ in $\Gamma$ and by hypothesis of case (II.2.2), $b\notin \Star(a)\cup \Star(c)\cup \Star(d) = \Link(a)\cup \Link(c)\cup \Link(d)$. We have that $\Link(x)\cap\Link(y)$ has $q$ elements, for $x\neq y$, $x,y\in \{a,c,d\}$. 
		As $\Gamma$ has no tetrahedrons, $\Link(\{a,b,c\})=\Link(a)\cap \Link(c)\cap \Link(d)$ is empty. Using the inclusion-exclusion principle,  we have $n-3l+3q$ choices for $b$.
		\item[(2.3)] $c\not \in \Star(\{a,d\}):$ we subdivide this case into the following disjoint subcases:
		\begin{itemize}
		    \item[(2.3.1)] $c=a$: This is impossible since $c\not\in \Star(\{a,d\})$.
		    \item[(2.3.2)] $c\in \Link(a)$: In this case, necessarily, $c \not \in \Star(d).$ Here we get $l-1-q$ choices for $c$. Also $b\notin \Link(a)\cup \Link(c)$ and we have $n-2l+q$ choices for $b$.
		    \item[(2.3.3)] $c\notin \Star(a):$ We do now again, three subcases:
		    \begin{itemize}
		        \item[(2.3.3.1)] $c=d$: This is impossible since $c\not\in \Star(\{a,d\})$.
		        \item[(2.3.3.2)] $c\in \Link(d)$: In this case $c\in \Link(d)\setminus\Star(a)$ here we get $l-1-q$ choices for $c$.  Also $b\notin \Link(a) \cup \Link(d)$ and we have $n-2l+q$ choices for $b$.
		        \item[(2.3.3.3)] $c\notin\Star(d):$ Here one has $|V\setminus(\Star(a)\cup\Star(d))|=n-2l+q$ choices for $c$.
		        Note that since $a,d$ span an edge and $b$ is not star of $a$ in (II.2.2) we have that $b$ can not be equal to $d$ neither to $a$. We subdivide this case into the following disjoint subcases:
			    \begin{enumerate}
			        \item[(2.3.3.3.1)] $b=d$: impossible.
				    \item[(2.3.3.3.2)] $b\in \Link(d)$: then $b\in \Link(d)\setminus\Star(a)$ and we have $l-1-q$ choices for $b$. 
				    \item[(2.3.3.3.3)] $b \not \in \Star(d):$  here $b\neq c$ to get a geodesic, and $b\notin \Star(a)\cup \Star(d)$. We have $n-2l+q-1$ choices for $b$.
			    \end{enumerate}
		    \end{itemize}
		\end{itemize}
	\end{itemize}
		Accounting for (2.1), (2.2) and (2.3), for a given vertex $a\in V$ we have:
		
	$$\sum_{\substack{c,b,d\in  \text{(II.2.2)}\\d\in \Link(a)}} E_{abcd}(z)=[l q (n-3l+3q)   +l(n-2l+q)(n+l-2q-3)]z^3 \mathcal{E}_v(z).$$
	\item[(3)]$d\not \in \Star(a)$: 
	\begin{itemize}
		\item[(3.1)] $b=d$: here we have $n-1-l$ choices for $d$, $1$ choice for $b$. As $abcb$ is a geodesic, $c$ does not belong to $\Star(b)$,  and we have  $n-l-1$ choices for $c$.
		\item[(3.2)] $b\neq d$: here we split into  these following disjoint cases:
		\begin{itemize}
			\item[(3.2.1)] $c=a$: here we get $1$ choice for $c$. Since $\{a,b\} \not\subseteq \Star(c)$ and $c = a$ we get $b\not \in \Star(a)$. One has that $b,d$ are any pair of different vertices of $V\setminus \Star(a)$, and they are distinct, thus we have $(n-l-1)(n-l-2)$ possibilities for  $b$ and $d$.
			\item[(3.2.2)] $c \in \Link(a)$: here we obtain $l$ choices for $c$. 	The hypothesis of case II.2, $\{a,b\}\not\subseteq\Star(c)$ implies that $b\not\in \Star(c)$. We consider  the following disjoint subcases:
			\begin{itemize}
			    \item[(3.2.2.1)] $d =c$: which is impossible since $abcd$ is geodesic.
				\item[(3.2.2.2)] $d \in \Link(c)$: here get $|\Link(c)\setminus \Star(a)|=l-1-q$ choices for $d$.
			 We have $|V\setminus (\Star(a)\cup\Star(c))|=n-2l+q$ for $b$.
				\item[(3.2.2.3)] $d \not \in \Star(c)$: we get $|V\setminus(\Star(a)\cup\Star(c))|=n-2l+q$ choices for $d$ and $|V\setminus(\Star(a)\cup \Star(c)\cup\{d\})|= n-2l+q-1$ choices for $b$.
			\end{itemize}
			\item[(3.2.3)] $c \not \in \Star(a)$: here for $b,c,d$ we get $(n-l-1)(n-l-2)(n-l-3)$ choices as $b,c,d$ can be any vertex outside of $\Star(a)$, $b\neq c$, $c\neq d$ because $abcd$ is geodesic, and $b\neq d$ by hypothesis.
		\end{itemize}
	\end{itemize}
			Accounting for (3.1) and (3.2), we obtain:
			
			$$\sum_{\substack{c,b,d\in  \text{(II.2.2)}\\d\not\in \Star(a)}} E_{abcd}(z)= [(n-1-l)^2+(n-l-1)(n-l-2)^2+l(n-2l+q)(n-l-2)] z^3\mathcal{E}_v(z).$$
\end{itemize}
Ultimately, in case (II.2.2), we obtain 

$$\sum_{\text{(II.2.2)}} E_{abcd}(z) = [(n-l-1)^3 + 2l(n-2l+q)(n-q-2)+lq(n-3l+3q)]z^3\mathcal{E}_v(z).$$

We finally collect all these cases together, and we conclude that:

\begin{align*}
\sum_{(a,b,c,d) \in (V\Gamma)^4} E_{abcd}(z) 
& = (q-1)z\mathcal{E}_\Delta (z)\\[-3pt]
& + (l-2)z\mathcal{E}_\Delta (z)\\[3pt]
& + [(l-1-q)^2 + q(l-2q)]z^2\mathcal{E}_e(z)\\[3pt]
& + (n-3)z\mathcal{E}_\Delta (z)\\[3pt]
& + [(l-1-q)(n-q-2)+q(l-2q)]z^2\mathcal{E}_e(z)\\[3pt]
& + [2(n-l-1) + q(n-q-l) + (n-q-2)(n-q-3)]z^2\mathcal{E}_e(z)\\[3pt]
& + [(n-l-1)^3 + 2l(n-2l+q)(n-q-2)+lq(n-3l+3q)]z^3\mathcal{E}_v(z)
\end{align*}

After grouping similar expressions we obtain:
\begin{align*}
\sum_{(a,b,c,d) \in (V\Gamma)^4} E_{abcd}(z) & =
(n+ q+l-6)z\mathcal{E}_\Delta (z)\\
& + (l^2 + ln + n^2 - 2lq - 2nq - 2q^2 - 6l - 4n + 10q + 7)z^2\mathcal{E}_e(z)\\
& + [(n-l-1)^3+2l(n-2l+q)(n-q-2)+lq(n-3l+3q)]z^3\mathcal{E}_v(z)
\end{align*} 
which expresses equation \eqref{E_4} of the main theorem.


\paragraph{Acknowledgments} 
The authors are grateful to the anonymous referee for carefully reading the first version of the manuscript and giving us several useful suggestions to improve the notation and presentation of the paper.

Yago  Antol\'{i}n  acknowledges  partial  support  from  the  Spanish  Government  through grants number MTM2017-82690-P, and through the ”Severo Ochoa Programme for Centres of Excellence in R\&D” (SEV-2015-0554) and (CEX2019-000904-S).

Islam Foniqi is a member of INdAM—GNSAGA, and gratefully acknowledges support from the Department of Mathematics of the University of Milano-Bicocca, and the Erasmus Traineeship grant 2020-1-IT02-KA103-078077.
\bibliography{main}
\bibliographystyle{plain}

\noindent\textit{\\ Yago Antol\'{i}n,\\
Fac. Matem\'{a}ticas, Universidad Complutense de Madrid and \\ 
Instituto de Ciencias Matem\'aticas, CSIC-UAM-UC3M-UCM\\
Madrid, Spain\\}
{email: yago.anpi@gmail.com}

\noindent\textit{\\ Islam Foniqi,\\
Università degli Studi di Milano - Bicocca\\ 
Milan, Italy\\}
{email: i.foniqi@campus.unimib.it}
\end{document}